\documentclass[11pt]{article}

\usepackage{lmodern} % evita warning
\usepackage{mathrsfs} % for script fonts
\usepackage{amsfonts}
\usepackage{amsthm,amsmath}
\usepackage{amsrefs}
\usepackage{color}
\usepackage{todonotes}
\newtheorem{theorem}{Theorem}[section]

\newtheorem{lemma}{Lemma}[section]
\newtheorem{remark}{Remark}[section]
\newtheorem{example}{Example}[section]

\newtheorem{assumption}{Assumption}

\usepackage{hyperref}

\usepackage{cleveref}
\MakeRobust{\Cref}

\begin{document}
\title{Large deviation principles and functional limit theorems in the deep
	limit of wide random neural networks\footnote{The authors acknowledge
	the support of MUR Excellence Department Project awarded to the Department
	of Mathematics, University of Rome Tor Vergata (CUP E83C23000330006),
    of University of Rome Tor Vergata (CUP E83C25000630005) Research
    Project METRO, and of INdAM-GNAMPA.}}
\author{Simmaco Di Lillo\thanks{Address: Dipartimento di Matematica,
	Università di Roma Tor Vergata, Via della Ricerca Scientifica,
	I-00133 Rome, Italy. e-mail: \texttt{dilillo@mat.uniroma2.it}}\and
	Claudio Macci\thanks{Address: Dipartimento di Matematica,
	Università di Roma Tor Vergata, Via della Ricerca Scientifica,
	I-00133 Rome, Italy. e-mail: \texttt{macci@mat.uniroma2.it}}\and
	Barbara Pacchiarotti\thanks{Address: Dipartimento di Matematica,
    Università di Roma Tor Vergata, Via della Ricerca Scientifica,
	I-00133 Rome, Italy. e-mail: \texttt{pacchiar@mat.uniroma2.it}}}
%\date{}
\maketitle
\begin{abstract}
	This paper studies large deviation principles and weak convergence, 
	both at the level of  finite-dimensional distributions  and in functional form, 
	for a class of continuous, isotropic, centered Gaussian random
	fields  defined on the unit sphere.    The covariance functions of these
	fields evolve recursively through a nonlinear map induced by an
	activation function, reflecting the statistical dynamics of infinitely
	wide random neural networks as depth increases. 
     We consider two types of centered fields, obtained by subtracting either the value at the
	North Pole or the spherical average.  According to the behavior of the derivative at $t=1$ of the associated covariance function, we identify three regimes: low disorder, sparse, and high disorder.
 In the low-disorder regime, we establish functional large deviation  principles and weak
	convergence results. In the sparse regime, we obtain large deviation
	principles and weak convergence for finite-dimensional distributions,
	while both properties fail at the functional level sense due to the emergence of discontinuities
	in the covariance recursion. 
\end{abstract}

\section{Introduction}
In recent years, the probabilistic analysis of neural networks has provided significant insights into
their asymptotic behaviour as the number of neurons per layer tends to infinity. In this limit, the outputs of deep feedforward random networks converge to Gaussian processes, whose  covariance functions evolve through a nonlinear recursion determined by the activation
function and initialization parameters. In this regime, often referred to as
kernel regimes or infinite-width limit, feedforward neural networks admit a probabilistic representation in terms of  Gaussian processes completely characterized by their associated positive definite kernels
\cite{RasmussenWilliams, Neal, Williams, DanielyFrostigSinger, Lee-et-al, deGMatthews-et-al, Hanin,
CammarotaMarinucciSalviVigogna, Balasubramanian-et-al, FavaroHaninMarinucciNourdinPeccati,
JacotGabrielHongler}.
Moreover, the kernel corresponding to a deep architecture admits a recursive
compositional structure: it can be written  as the $L$-fold iteration of a fixed base kernel, where
$L$ denotes the network depth. The connection between Gaussian random fields and deep neural architectures
has motivated several studies aimed at understanding the statistical and geometric properties of neural
representations in high dimensions. 

A natural question in this context is to understand the asymptotic fluctuations of these random fields as
the depth increases. In this direction,  recent works have investigated  large deviation principles
for models of interacting neurons and random neural networks in mean-field regimes, establishing large deviation principles for empirical processes or recursive Gaussian systems and  providing a rigorous
probabilistic framework to quantify rare fluctuations in high-dimensional neural models (see, for instance,
\cite{HirschWillhalm, MacciPacchiarottiTorrisi, AndreisBassettiHirsch, Vogel}). These developments motivate
our present investigation.

It is known (see, e.g., \cite{DililloMarinucciSalviVigogna}) that the qualitative behaviour of  deep random neural networks is governed by the properties of the iterated covariance function. In particular, depending on the value of the derivative $\kappa^\prime(1)$, one can
distinguish the regimes, referred as the \emph{low-disorder}, \emph{sparse}, and \emph{high-disorder} phases.
A key motivation for our analysis comes from Theorem 3.3 in \cite{DililloMarinucciSalviVigogna}, where  convergence to zero of the associated Gaussian fields is established in then low-disorder and sparse regimes after centering at the North Pole. In the present work, we aim to refine these results within a large deviation framework and consider, in addiction, a second centering obtained by subtracting the spherical average of the field.
%In that theorem one has these Gaussian fields after subtracting the value at the North Pole (\emph{Case 1}); here we also consider these Gaussian fields after subtracting its integral over the sphere (\emph{Case 2}).

In this paper, we establish large deviation principles  for finite-dimensional distributions in both the low-disorder and sparse regimes, for each of the two centering procedures considered. We  also prove some finite-dimensional weak convergence results. A fundamental distinction arises between the two regimes. In the low-disorder case, we obtain functional large deviation principles and functional weak convergence, whereas in the sparse
case both fail  at the functional level due to the emergence of discontinuities in the covariance functions. We also present a study for the high-disorder
regime.  The difference between low-disorder and sparse scenarios can be attributed to the relative weight of the network
initialization as the depth increases. In the low-disorder regime, the Gaussian field exhibits an exponential decay toward
 triviality, both at the level of finite-dimensional distributions and in functional spaces. In contrast, in the sparse regime the decay us polynomial, and although pointwise
convergence to zero still holds, the initialization weight prevents the validity of functional convergence results. This divergence reflects a deeper structural contrast:  depth erases information uniformly in the low-disorder regime, its effect remains visible in the sparse case.

The structure of the paper is as follows. In \Cref{sec:preliminaries} we recall the main probabilistic and analytic preliminaries. Finite-dimensional and functional
results are presented in \Cref{sec:fdd} and in \Cref{sec:functional}, respectively. The high-disorder regime is discussed in  \Cref{sec:high-disorder}. For the sake of readability, all
proofs are deferred to Appendix \ref{app}.

Overall, this paper contributes to extending classical asymptotic theory for Gaussian random fields to
recursive and nonlinear settings, and clarify the probabilistic mechanisms that govern how depth influences the
behaviour of wide random neural networks. By linking spectral observations
with probabilistic limits, we provide a framework for quantifying the degeneration of signal propagation
in very deep architectures and for understanding  how depth and initialization regimes affect
the expressive and statistical properties of random neural networks.  This framework also set the stage for future developments involving higher-order fields, such as those incorporating derivatives \cite{DililloMarinucciSalviVigogna}), as well as extensions to non-Gaussian models.

\section{Preliminaries}\label{sec:preliminaries}
In this section we present some preliminaries on large deviations, on isotropic random fields on the sphere,
and on random neural networks.

\subsection{Large deviations}
We refer to \cite{DemboZeitouni} as a reference on this topic. We start with the definition of large
deviation principle (LDP from now on). We refer to a sequence of random variables
$\{X_L\}_{L\geq 1}$ (taking values on some topological space $\mathcal{X}$, and defined on the same probability
space $(\Omega,\mathcal{F},P)$), letting $L\to\infty$, because this formulation is more direct and better
suited for what follows.
A lower semi-continuous function $I:\mathcal{X}\to[0,\infty]$ is called rate function, and it is said to be
good if all its level sets $\{\{x\in\mathcal{X}:I(x)\leq\eta\}:\eta\geq 0\}$ are compact. Then
$\{X_L\}_{L\geq 1}$ satisfies the LDP with speed $v_L\to\infty$ and rate function $I$ if
$$\limsup_{L\to\infty}\frac{1}{v_L}\log P(X_L\in C)\leq-\inf_{x\in C}I(x)\ \mbox{for all closed sets}\ C\subset\mathcal{X}$$
and
$$\liminf_{L\to\infty}\frac{1}{v_L}\log P(X_L\in O)\geq-\inf_{x\in O}I(x)\ \mbox{for all open sets}\ O\subset\mathcal{X}.$$
In general, we deal with centered Gaussian sequences. We shall prove large deviation results for finite dimensional
distributions and we consider some applications of the G\"artner Ellis Theorem (see, e.g., Theorem 2.3.6 in
\cite{DemboZeitouni}). We also present functional results for sequences of continuous random fields defined on the
sphere. Here,  the rate functions will be expressed
in terms of a suitable supremum of rate functions for finite dimensional distributions (see, e.g., \cite{Baldi});
this will be clear looking at the statement of Theorem \ref{prop:low-disorder-functional}.
Finally, we also refer to the \emph{contraction principle} (see, e.g., Theorem 4.2.1 in \cite{DemboZeitouni}) and the Dawson G\"artner Theorem (Theorem 4.6.1 in \cite{DemboZeitouni}).

\subsection{Isotropic random fields on the sphere}

Throughout this paper we use the symbol $C(\mathbb{S}^d)$ to denote the space of all real-valued continuous
functions on $\mathbb{S}^d$. We recall some well known properties of continuous, isotropic, centered Gaussian random
fields $(T(x))_{x\in\mathbb{S}^d}$ defined on a probability space $(\Omega,\mathcal{F},P)$; see, for instance,
\cites{Yadrenko1983,MarinucciPeccati}.
We denote by $\lambda$ the Lebesgue measure on $\mathbb{S}^d$, by $\omega_d=\lambda(\mathbb{S}^d)$ its total mass, and we make use of the following spectral representation
\begin{equation}\label{eq:isotropic-random-field}
	T(x,\omega)
	= a_{00}(\omega) Y_{00}(x)
	+ \sum_{\ell=1}^{\infty} \sum_{m=1}^{n_{\ell,d}} a_{\ell m}(\omega)\, Y_{\ell m}(x),
\end{equation}
which holds in $L^{2}(\Omega\times\mathbb{S}^d)$, that is,
\[
\lim_{\ell_*\to\infty}
\int_{\Omega}\int_{\mathbb{S}^d}
\left|
T(x,\omega) -
\left(
a_{00}(\omega)Y_{00}(x)
+ \sum_{\ell=1}^{\ell_*}\sum_{m=1}^{n_{\ell,d}} a_{\ell m}(\omega) Y_{\ell m}(x)
\right)
\right|^2
\lambda(dx)\,P(d\omega)
=0.
\]

Here:
\begin{itemize}
	\item $n_{\ell,d}=\dfrac{2\ell+d-1}{\ell}\binom{\ell+d-2}{\ell-1}$ for $\ell\ge 1$;
	\item $\mathcal{A}=\{(0,0)\}\cup\{(\ell,m):\ell\ge1,\ m=1,\dots,n_{\ell,d}\}$;
	\item $\{a_{\ell m}\}_{(\ell,m)\in\mathcal{A}}$ are independent centered Gaussian random variables with
	$\mathbb{E}[a_{\ell m}^2]=C_\ell\ge0$;
	\item $\{Y_{\ell m}\}_{(\ell,m)\in\mathcal{A}}$ is an $L^2(\lambda)$ orthonormal basis of real--valued spherical harmonics,
	each $Y_{\ell m}$ being an eigenfunction of the Laplace-Beltrami operator with eigenvalue $-\ell(\ell+d-1)$.
\end{itemize}

From the spectral representation and orthonormality we immediately obtain
\[
\frac{1}{\omega_d}\int_{\mathbb{S}^d} T(x)\, \lambda(dx)
= \frac{a_{00}}{\sqrt{\omega_d}},
\]
since $Y_{00}(x)=Y_{00}\equiv \omega_d^{-1/2}$.

Moreover, isotropy implies the existence of a continuous function $\kappa:[-1,1]\to\mathbb{R}$ such that
\begin{equation}\label{eq:def-kappa}
	\mathrm{Cov}(T(x),T(y)) = \kappa(\langle x,y\rangle)
	\qquad \text{for all } x,y\in\mathbb{S}^d,
\end{equation}
and in particular $\mathrm{Var}[T(x)]=\kappa(1)$. Using the addition formula for spherical harmonics~\cite{AtkisonHan}, the covariance function admits the expansion
\[
\kappa(t)
= \sum_{\ell\ge0} \frac{C_\ell n_{\ell,d}}{\omega_d}\, G_{\ell,d}(t),
\]
where $G_{\ell,d}$ is the normalized Gegenbauer polynomial with  degree $\ell$,  uniquely determined by the value
$G_{\ell,d}(1)=1$ and the orthogonality relation (see \cite{Szego1975})
\[
\int_{-1}^1 G_{\ell,d}(t)\, G_{\ell',d}(t)\, (1-t^2)^{\frac{d}{2}-1} \, dt = 0,
\qquad \ell\neq\ell'.
\]
For convenience, we set
\[
D_\ell = \frac{C_\ell n_{\ell,d}}{\omega_d},
\]
so that the covariance satisfies
\[
\kappa(1) = \sum_{\ell\ge0} D_\ell.
\]

\subsection{Random neural networks}

We recall some preliminary notions concerning random neural networks; see, for instance,
Section~2.2 and Proposition~A.6 in~\cite{DililloMarinucciSalviVigogna}.
A random neural network is defined recursively as follows.
For the input layer, we set
\[
Z_{0,i}(x) = \sum_{j=1}^{m_0} W^{(0)}_{ij} x_j + b^{(1)}_i,
\qquad i=1,\dots,m_1,
\]
and, for $s=1,\dots,L$,
\[
Z_{s,i}(x) =
\sum_{j=1}^{m_s} W^{(s)}_{ij}\,\sigma\!\big(Z_{s-1,j}(x)\big)
+ b^{(s+1)}_i,
\qquad i=1,\dots,m_{s+1}.
\]

The random variables appearing above are independent and Gaussian with suitable centered laws:
the biases $b^{(s)}_i$, $s=1,\dots,L$, have variance $\Gamma_b \ge 0$,
the initial weights $W^{(0)}_{ij}$ have variance $\Gamma_{W_0}>0$,
and the weights $W^{(s)}_{ij}$ for $s\ge1$ have variance $\Gamma_W/m_s >0$.
We denote by $\sigma:\mathbb{R}\to\mathbb{R}$ the activation function;
$L$ is the total number of hidden layers;
$m_0=d+1$ is the input dimension, $m_{L+1}$ is the output dimension, and
$m_1,\dots,m_L$ are the widths of the hidden layers.

\medskip

The \emph{wide-limit} of the random neural network, i.e. the limit as $m_1,\dots,m_L \to \infty$, has been studied in many works; see, e.g., \cites{Neal, Hanin, CammarotaMarinucciSalviVigogna, Balasubramanian-et-al, FavaroHaninMarinucciNourdinPeccati} and the references  therein.
More precisely, one shows that the family
\[
\big\{ (Z_{L,i}(x))_{x\in \mathbb S^d} \big\}_{i=1,\dots,m_{L+1}}
\]
converges weakly, as the widths diverge, to a Gaussian random vector with i.i.d.\ components
\[
\big\{ (Z^*_{L,i}(x))_{x\in\mathbb{S}^d} \big\},
\]
whose covariance function can be defined recursively with respect to the depth~$L$.
In particular, since the limiting components are i.i.d., it suffices to study a single univariate limiting component, which we denote by $(T_L(x))_{x\in \mathbb S^d}$.

Under the standard calibration conditions,
	\[
	\mathbb{E}[\sigma(Z)^2] < \infty,
	\qquad
	\Gamma_b + \Gamma_{W_0} =1,
	\qquad
	\Gamma_W\,\mathbb{E}[\sigma(Z)^2] + \Gamma_b = 1,
	\]
	for a standard normal random variable $Z$, the limiting covariance function is given by
	\[
	\mathrm{Cov}\big(T_L(x),T_L(y)\big)
	= \kappa_L\!\left(\Gamma_b + \Gamma_{W_0}\,\langle x,y\rangle \right),
	\]
	where  $\kappa_L$ denotes the $L$-fold composition of $\kappa$, i.e., 
    \[
\kappa_L = \underbrace{\kappa\circ\cdots\circ\kappa}_{L\ \text{times}}
\]
and
	\[
	\kappa(t)
	= \Gamma_W\,\mathbb{E}\!\left[\sigma(\widehat Z_1)
	\,\sigma\!\big(t\,\widehat Z_1 + \sqrt{1-t^2}\,\widehat Z_2\big)\right]
	+ \Gamma_b, \qquad t \in [-1,1]
	\]
	with $\widehat Z_1,\widehat Z_2$ i.i.d.\ standard normals.
	See~\cite{DililloMarinucciSalviVigogna} for more details.
If $\sigma$ belongs to the Malliavin--Sobolev space, it admits an $L^2$ Hermite expansion, i.e.
	\[
	\sigma(x)
	= \sum_{q\ge0} \frac{\widehat J_q(\sigma)}{q!}\, H_q(x),
	\qquad
	\sum_{q\ge1} \frac{|\widehat J_q(\sigma)|}{(q-1)!} < \infty.
	\]
Using the Diagram Formula (see \cite{MarinucciPeccati}*{Prop.~4.15}), one obtains
	\begin{equation}\label{eq:kappa-for-sigma}
		\kappa(t)
		= \sum_{q\ge0} \frac{(J_q(\sigma))^2}{q!}\, t^q,  \quad J_q(\sigma) = \begin{cases} \sqrt{\Gamma_W \widehat J_0(\sigma)^2 + \Gamma_b} & \text{ if } q=0 ,\\
        \sqrt{\Gamma_W} \widehat{J_q}(\sigma) & \text{ otherwise}.
        \end{cases}
	\end{equation}
We further assume that $\sigma$ is non-polynomial; otherwise each layer would yield a polynomial transformation of the input and their composition would remain polynomial, preventing universal approximation (see \cites{Cybenko, HornikStinchcombeWhite}).
Equivalently, $\widehat J_q(\sigma)\neq 0$ for infinitely many values of $q$.
		The above summability condition ensures that $\kappa\in C^1([-1,1])$.
        
From now on, for the sake of clarity, 
we restrict to the case $\Gamma_b = 0$.
Then $\Gamma_{W_0}=1$ and
\[\Gamma_W = \mathbb{E}[\sigma(Z)^2]^{-1}
\]
Although a version of the results in this work could also be developed for the general case $\Gamma_b\in(0,1)$,
the associated notation and arguments become substantially more involved.

\subsection{Three regimes of spectral complexity}

Motivated by the structure of the Gaussian fields that arise as infinite-width limits
of random neural networks, we shall focus in this work on a class of sequences of
continuous, isotropic, centered Gaussian random fields
\(
\big\{(T_L(x))_{x\in \mathbb S^d} \big\}_{L\ge 1},
\)
where, for every $L\ge 1$,
\begin{equation}\label{KL}
	\mathrm{Cov}(T_L(x),T_L(y))=\kappa_L(\langle x,y\rangle),
	\qquad x,y\in\mathbb{S}^d.
\end{equation}
Here $\kappa:[-1,1]\to\mathbb{R}$ is a $C^1([-1,1])$ function that admits the power-series
	representation in \eqref{eq:kappa-for-sigma}
	%\begin{equation}\label{eq:power-series}
	%	\kappa(u)=\sum_{q\ge 0} a_q^2\, u^q,
	%\end{equation}
	with $J_q(\sigma)\neq 0$ for infinitely many $q$, and normalized so that $$\kappa(1)=1.$$
Moreover, $\kappa_L$ again denotes the $L$-fold composition of~$\kappa$.
\medskip

The asymptotic behaviour of the fields
\(
\big\{(T_L(x))_{x\in\mathbb{S}^d}\big\}_{L\ge 1}
\)
as \(L\to\infty\) is investigated in
\cite{DililloMarinucciSalviVigogna}, where three regimes are identified, depending
solely on the value of the derivative of the covariance kernel at~1:
\begin{itemize}
	\item \emph{low-disorder}, if $\kappa'(1)<1$;
	\item \emph{sparse}, if $\kappa'(1)=1$;
	\item \emph{high-disorder}, if $\kappa'(1)>1$.
\end{itemize}

\begin{remark}\label{rem:k'(1)-role}
	The quantity $\kappa'(1)$ determines the local behaviour of the covariance function
	for nearby points and thus controls the rate of decay of correlations.
	It is also directly related to the gradient variance; in particular, it governs the
	local smoothness of the field.
	Furthermore, in Kac--Rice-type formulae, $\kappa'(1)$ plays a crucial role in the
	asymptotics of the expected number of critical points or nodal sets.
	Hence, $\kappa'(1)$ acts as a bridge between the spectral properties of the random
	field and its geometric or topological features.
	These aspects are discussed in detail in
	\cites{Dilillo, DililloMarinucciSalviVigogna-fractal}; see also the classical
	references \cites{AdlerTaylor,AzaisWschebor}.
\end{remark}

\medskip

Our aim is to study, from the viewpoint of large deviations, certain convergences
to zero that have been established in the literature
(see Lemma~\ref{lem:statement}).
More precisely, we regard the sequences appearing in Lemma~\ref{lem:statement} as
random fields, that is, with \(x\) ranging over~$\mathbb{S}^d$.
We begin by considering an arbitrary finite
set of   points
$x_1,\ldots,x_m\in\mathbb{S}^d$ with $x_i\neq x_j$ for $i\neq j$, where $m\ge 1$.
In this setting, we deal with sequences taking values in $\mathbb{R}^m$ that
converge to the null vector~$\underline{0}$, or, from a functional point of view, with
sequences of fields converging to the identically zero function on~$\mathbb{S}^d$.

\begin{lemma}\label{lem:statement}
	If $\kappa^\prime(1)\le 1$, then the following sequences converge to zero in $L^2(\Omega)$, as $L\to\infty$, for each fixed
	$x\in\mathbb{S}^d$:
	\begin{enumerate}
		\item $\{T_L(x)-T_L(N)\}_{L\ge 1}$, where $N$ denotes the North Pole of
		$\mathbb{S}^d$;
		\item
		$\displaystyle
		\left\{
		T_L(x)-\frac{1}{\omega_d}\int_{\mathbb{S}^d}T_L(y)\,\lambda(dy)
		\right\}_{L\ge 1}
		= \left\{T_L(x)-\frac{a_{00,L}}{\sqrt{\omega_d}}\right\}_{L\ge 1}.
		$
	\end{enumerate}
\end{lemma}

\begin{proof}
	The convergence in item~1 follows from Theorem~3.3 and identity~(3.7)
	in~\cite{DililloMarinucciSalviVigogna}.
	The convergence in item~2 is a consequence of the results summarized in Table~1 of
	\cite{DililloMarinucciSalviVigogna}.
\end{proof}

\section{Results for finite-dimensional distributions}\label{sec:fdd}

In this section we consider the sequences introduced in \Cref{lem:statement} and
study their finite-dimensional distributions.

Throughout, we let
\[
\underline{x}=(x_1,\ldots,x_m)\in (\mathbb{S}^d)^m ,
\]
denote configurations of pairwise distinct points on the sphere as previously specified.
We also write
\[
T_L(\underline{x})=\big(T_L(x_1),\ldots,T_L(x_m)\big)\in\mathbb{R}^m,
\]
and denote by \(\boldsymbol{1}\in\mathbb{R}^m\) the vector whose coordinates are all equal to~1.

\medskip

A standing assumption in this work concerns the regularity of the covariance function~\(\kappa\).
Before stating our main results, we formulate this condition explicitly.

\begin{assumption}[Covariance Regularity]\label{ass:A}
	The covariance function \(\kappa:[-1,1]\to\mathbb{R}\) satisfies, as \(t\to1^{-}\),
	\begin{equation}\label{eq:sort-of-Taylor-formula}
		1 - \kappa(t)
		= \kappa'(1)(1-t)
		\;-\;
		c (1-t)^{\rho}
		\;+\;
		o\!\big((1-t)^{\rho}\big),
	\end{equation}
	for some \(c\in\mathbb{R}\setminus\{0\}\) and some \(\rho>1\).
\end{assumption}

\begin{remark}
	Assumptions of this type appear in several references; see, for instance,
	\cite{BiettiBach} (Theorem~1 and the explicit formulas for ReLU networks in Section~3.2),
	and \cite{DililloMarinucciSalviVigogna-fractal} (Definition~3.1, where the exponent \(\rho\)
	is referred to as the \emph{Covariance Regularity Index}).
	In the latter reference it is further shown that, whenever \(\rho>1\), the associated limiting
	fields are almost surely of class \(C^1(\mathbb S^d)\).

	We stress that in the low-disorder regime our LDP for finite-dimensional distributions
	does not require Assumption~\ref{ass:A}.
	In this case, however, we are unable to exclude the possibility that the corresponding
	rate function may be trivial (see \Cref{rem:trivial-cases-revised} for details). By contrast, in the sparse regime Assumption~\ref{ass:A} becomes crucial.
\end{remark}

We are now ready to formulate the large deviation principles.
Each of the two theorems below establishes two distinct LDPs, corresponding to the
sparse and the low-disorder regimes.
In both theorems, two limiting functions appear, one for each regime.
Their existence is guaranteed by
\Cref{lem:low-disorder,lem:sparse}
which concern respectively the low-disorder and the sparse regimes.

\begin{theorem}\label{teo:low-disorder-1}
	Let $\kappa'(1)\le 1$; moreover if $\kappa'(1)=1$ we assume also Assumption~\ref{ass:A}.
	Then the sequence
	\[
	\big\{\,T_L(\underline{x}) - T_L(N)\,\boldsymbol{1}\,\big\}_{L\ge 1}
	\]
	satisfies an LDP with speed $v_L$ and good rate function
	$I_{1;\mathfrak{g}}(\cdot;\underline{x})$, given by
	\begin{equation}\label{eq:rf-1-low-disorder}
		I_{1;\mathfrak{g}}(\underline{y};\underline{x})
		=
		\sup_{\underline{\theta}\in\mathbb{R}^m}
		\left\{
		\langle \underline{\theta},\underline{y}\rangle
		-\frac12 \langle \underline{\theta},
		B_{1;\mathfrak{g}}(\underline{x})\,\underline{\theta} \rangle
		\right\}, \quad \underline{y}\in \mathbb{R}^m,
	\end{equation}
	where $\mathfrak{g}:[-1,1]\to\mathbb{R}$,  $B_{1;\mathfrak{g}}(\underline{x})$ is the matrix
	\begin{equation}\label{eq:matrixB1-low-disorder}
		B_{1;\mathfrak{g}}(\underline{x})
		=
		\big(
		\mathfrak{g}(\langle x_i,N\rangle)
		+\mathfrak{g}(\langle x_j,N\rangle)
		-\mathfrak{g}(\langle x_i,x_j\rangle)
		\big)_{i,j=1}^m ,
	\end{equation}
and:
	\begin{itemize}
		\item if $\kappa'(1)<1$, then $v_L = (\kappa'(1))^{-L}$ and $\mathfrak{g}=\mathfrak{L}$, where
		\[
		\mathfrak{L}(t)
		:=\lim_{L\to\infty} (\kappa'(1))^{-L}\bigl(1-\kappa_L(t)\bigr);
		\]

		\item if $\kappa'(1)=1$, then $v_L = L^{1/(\rho-1)}$ and $\mathfrak{g}=\mathfrak{S}$, where
		\[
		\mathfrak{S}(t)
		:=\lim_{L\to\infty} L^{1/(\rho-1)}\bigl(1-\kappa_L(t)\bigr).
		\]
	\end{itemize}
\end{theorem}

\begin{theorem}\label{teo:low-disorder-2}
	Under the same assumptions as in \Cref{teo:low-disorder-1}, the sequence
	\[
	\left\{
	T_L(\underline{x})
	-  \frac{a_{00,L}}{\sqrt{\omega_d}}
\boldsymbol{1}
	\right\}_{L\ge 1}
	\]
	satisfies an LDP with speed $v_L$ and good rate function
	$I_{2;\mathfrak{g}}(\cdot;\underline{x})$, given by
	\begin{equation}\label{eq:rf-2-low-disorder}
		I_{2;\mathfrak{g}}(\underline{y};\underline{x})
		=
		\sup_{\underline{\theta}\in\mathbb{R}^m}
		\left\{
		\langle\underline{\theta},\underline{y}\rangle
		-\frac12
		\langle\underline{\theta},B_{2;\mathfrak{g}}(\underline{x})\underline{\theta}\rangle
		\right\}, \quad \underline{y}\in \mathbb{R}^m,
	\end{equation}
	where $B_{2;\mathfrak{g}}(\underline{x})$ is the matrix
	\begin{equation}\label{eq:matrixB2-low-disorder}
		B_{2;\mathfrak{g}}(\underline{x})
		=
		\left(
		\frac{\displaystyle\int_{-1}^1 \mathfrak{g}(t)\,(1-t^2)^{\frac d2 -1}\,dt}
		{\displaystyle\int_{-1}^1 (1-t^2)^{\frac d2 -1}\,dt}
		-
		\mathfrak{g}(\langle x_i,x_j\rangle)
		\right)_{i,j=1}^m ,
	\end{equation}
	and $v_L$ and $\mathfrak{g}$ are as in \Cref{teo:low-disorder-1}.
\end{theorem}

\Cref{lem:low-disorder,lem:sparse} show that the functions $\mathfrak{L}$ and $\mathfrak{S}$ are well defined.
We also show that the limiting function vanishes only on the set of points mapped to~1 by $\kappa$, namely
\[
\mathcal{I}(\kappa):=\{t\in[-1,1]:\kappa(t)=1\}.
\]
The following lemma characterizes this set completely.

\begin{lemma}\label{lem:set-of-attractive-points-of-iterations}
	Assume that $\kappa'(1)\le 1$. Then exactly two situations may occur: either
	$\mathcal{I}(\kappa)=\{1\}$ or $\mathcal{I}(\kappa)=\{-1,1\}$.
	Moreover, $\mathcal{I}(\kappa)=\{-1,1\}$ if and only if $\kappa$ is symmetric, i.e.
	$\kappa(t)=\kappa(-t)$ for all $t\in[-1,1]$.
\end{lemma}

We are now in a position to state the existence of the two limiting functions.

\begin{lemma}\label{lem:low-disorder}
	Assume that $\kappa'(1)<1$. Then the limit
	\[
	\mathfrak{L}(t)
	=\lim_{L\to\infty} (\kappa'(1))^{-L}\bigl(1-\kappa_L(t)\bigr)
	\]
	exists in $[0,\infty)$.
	Moreover, under Assumption~\ref{ass:A}, we have $\mathfrak{L}(t)=0$ if and only if $t\in\mathcal{I}(\kappa)$.
\end{lemma}

\begin{lemma}\label{lem:sparse}
	Assume that $\kappa'(1)=1$ and Assumption~\ref{ass:A} holds. Then
	\[
	\mathfrak{S}(t)
	=
	\lim_{L\to\infty} L^{1/(\rho-1)}\bigl(1-\kappa_L(t)\bigr)
	=
	\begin{cases}
		0, & t\in\mathcal{I}(\kappa),\\[4pt]
		h_{c,\rho}, & t\notin\mathcal{I}(\kappa),
	\end{cases}
	\]
	where $h_{c,\rho}=(c(\rho-1))^{-1/(\rho-1)}$ and $c>0$ is the coefficient of $(1-t)^\rho$ in
	\eqref{eq:sort-of-Taylor-formula}.
\end{lemma}

\begin{remark}[Assumption~\ref{ass:A} avoids triviality in the low-disorder case]
	\label{rem:trivial-cases-revised}
	Assume that $\kappa'(1)<1$. We observe that the rate functions
	$I_{1,\mathfrak{L}}$ and $I_{2,\mathfrak{L}}$ vanish at
	$\underline{y}=\underline{0}\in\mathbb{R}^m$, and therefore one cannot exclude
	a priori the possibility that the resulting rate function $I$ is trivial, namely
	\[
	I(\underline{y})=
	\begin{cases}
		0, & \underline{y}=\underline{0},\\[4pt]
		\infty, & \text{otherwise}.
	\end{cases}
	\]
	This scenario cannot occur if the matrices $B_{1,\mathfrak{L}}$ and
	$B_{2,\mathfrak{L}}$ have non-zero rank.
	We now show that, under Assumption~\eqref{eq:sort-of-Taylor-formula}, both matrices, for $m>2$,
	necessarily contain non-zero entries.

	\begin{enumerate}
		\item The matrix $B_{1,\mathfrak{L}}$ has diagonal entries
		$2\mathfrak{L}(\langle x_i,N\rangle)$ for $i=1,\ldots,m$.
		If $x_i\notin\{ S,N\}$, where $S$ denotes the South Pole, then $\langle x_i,N\rangle\neq \pm1$, and by the second part
		of \Cref{lem:low-disorder} we have $\mathfrak{L}(\langle x_i,N\rangle)>0$.
		Thus, the corresponding diagonal entry of $B_{1,\mathfrak{L}}$ is non-zero.

		\item The matrix $B_{2,\mathfrak{L}}$ has diagonal entries equal to
		\[
		\frac{\displaystyle\int_{-1}^1
			\mathfrak{L}(t)(1-t^2)^{\frac{d}{2}-1}\,dt}
		{\displaystyle\int_{-1}^1 (1-t^2)^{\frac{d}{2}-1}\,dt}.
		\]
		Since $\mathfrak{L}(t)\ge 0$ for all $t$ and $\mathfrak{L} \neq 0$ in $(-1,1)$,
		this quantity is strictly positive, and therefore each diagonal entry is non-zero.
	\end{enumerate}
\end{remark}

Since in the sparse case the limiting function $\mathfrak{S}$ takes a particularly simple form, explicit expressions for the two rate functions $I_{1;\mathfrak{S}}$ and $I_{2;\mathfrak{S}}$ can be derived.
These formulas will also be instrumental in showing that, in the sparse regime, the two theorems above do not admit a functional counterpart (see~\Cref{prop:sparse-functional}).

\begin{remark}[Explicit formulation for $I_{1;\mathfrak{S}}$]\label{rem:rf-sparse-1}
	We distinguish two cases.

	\begin{itemize}
		\item Suppose first that $\kappa$ is not symmetric, so that
		$\mathcal{I}(\kappa)=\{1\}$.
		If none of the points $x_1,\ldots,x_m$ coincides with the North Pole~$N$,
		the matrix $B_{1;\mathfrak{S}}(\underline{x})=(b_{ij}(\underline{x}))_{i,j=1}^m$ is given by 
\[
	b_{ij}(\underline{x})
	=
	\begin{cases}
		2 h_{c,\rho}, & i=j,\\[4pt]
		h_{c,\rho}, & i\neq j,\ 
	\end{cases}
	\]
        %$$ (B_{1;\mathfrak{S}}(\underline{x}))_{ij} =  h_{c,\rho} (1 + \delta_{ij}), $$
        and hence
is invertible. One  can check that for $\underline{y}=(y_1,\ldots, y_m)$ we have
		\begin{equation}\label{eq:rf-1-sparse}
			I_{1;\mathfrak{S}}(\underline{y};\underline{x})
			=\frac{1}{2(m+1)h_{c,\rho}}
			\left(
			m\sum_{i=1}^m y_i^2
			-2\sum_{1\le i<j\le m} y_i y_j
			\right).
		\end{equation}

		If instead $x_i=N$ for some $i\in\{1,\ldots,m\}$ (we assume $m\ge2$ to avoid trivialities),
		let $\underline{y}_i$ and $\underline{x}_i$ denote the vectors obtained by removing the
		$i$-th components from $\underline{y}$ and $\underline{x}$, respectively.
		In this case the problem reduces to $m-1$ points, and
		\[
		I_{1;\mathfrak{S}}(\underline{y};\underline{x})
		=
		\begin{cases}
			I_{1;\mathfrak{S}}(\underline{y}_i;\underline{x}_i), & y_i=0,\\[4pt]
			\infty, & \text{otherwise}.
		\end{cases}
		\]

		\item Suppose now that $\kappa$ is symmetric, and hence
		$\mathcal{I}(\kappa)=\{-1,1\}$.
		Formula \eqref{eq:rf-1-sparse} remains valid whenever no two points in
		$\underline{x}$ are antipodal.
		If $x_i=N$ or $x_i=S$, and its antipodal
		point is not present among $x_1,\ldots,x_m$, then the same reduction to
		$m-1$ points as above applies.

		If, on the other hand, there are $2v$ antipodal points, say
		\[
		x_{i_1}=-x_{j_1},\;\ldots,\; x_{i_v}=-x_{j_v},
		\]
		let $\underline{y}_{\underline{i}_v}$ and $\underline{x}_{\underline{i}_v}$
		denote the vectors obtained by removing the components with indices
		$i_1,\ldots,i_v$.
		Then the problem reduces to a system with $m-v$ points, and
		\[
		I_{1;\mathfrak{S}}(\underline{y};\underline{x})
		=
		\begin{cases}
			I_{1;\mathfrak{S}}(\underline{y}_{\underline{i}_v};
			\underline{x}_{\underline{i}_v}),
			& y_{i_1}=y_{j_1},\ldots,y_{i_v}=y_{j_v},\\[4pt]
			\infty, & \text{otherwise}.
		\end{cases}
		\]

		Finally, if for some $h\in\{1,\ldots,v\}$ we have
		$(x_{i_h},x_{j_h})=(N,S)$ or $(S,N)$, then the reduction is to $m-v-1$ points:
		one must impose $y_{i_h}=y_{j_h}=0$; otherwise,
		$I_{1;\mathfrak{S}}(\underline{y}_{\underline{i}_v};
		\underline{x}_{\underline{i}_v})=\infty$.
	\end{itemize}
\end{remark}

\begin{remark}[Explicit formulation for $I_{2;\mathfrak{S}}$]\label{rem:rf-sparse-2}
	The matrix $B_{2,\mathfrak{S}}(\underline{x})=(b_{ij}(\underline{x}))_{i,j=1}^m$ appearing in
	\eqref{eq:matrixB2-low-disorder} takes the form
	\[
	b_{ij}(\underline{x})
	=
	\begin{cases}
		h_{c,\rho}, & i=j,\\[4pt]
		0, & i\neq j\ \text{and}\ \langle x_i,x_j\rangle\notin\mathcal{I}(\kappa),\\[4pt]
		h_{c,\rho}, & i\neq j\ \text{and}\ \langle x_i,x_j\rangle\in\mathcal{I}(\kappa).
	\end{cases}
	\]
	As before, we distinguish two cases according to Lemma~\ref{lem:set-of-attractive-points-of-iterations}.

	\begin{itemize}
		\item Suppose first that $\kappa$ is not symmetric, so that
		$\mathcal{I}(\kappa)=\{1\}$.
		In this case $B_{2,\mathfrak{S}}(\underline{x})$ is diagonal and invertible, and hence
		\begin{equation}\label{eq:rf-2-sparse}
			I_{2;\mathfrak{S}}(\underline{y};\underline{x})
			=\frac{1}{2h_{c,\rho}}\|\underline{y}\|^2.
		\end{equation}

		\item Now assume that $\kappa$ is symmetric, and therefore
		$\mathcal{I}(\kappa)=\{-1,1\}$.
		Formula \eqref{eq:rf-2-sparse} remains valid whenever $\underline{x}$ contains no
		antipodal pairs.
		If instead there are $2v$ antipodal points, say
		\[
		x_{i_1}=-x_{j_1},\;\ldots,\; x_{i_v}=-x_{j_v},
		\]
		let $\underline{y}_{\underline{i}_v}$ and $\underline{x}_{\underline{i}_v}$ denote the
		vectors obtained by removing the components with indices $i_1,\ldots,i_v$.
		In this case the problem reduces to one with $m-v$ points, and
		\[
		I_{2;\mathfrak{S}}(\underline{y};\underline{x})
		=
		\begin{cases}
			\displaystyle
			\frac{1}{2h_{c,\rho}}
			\left\|\underline{y}_{\underline{i}_v}\right\|^2,
			& y_{i_1}=y_{j_1},\;\ldots,\; y_{i_v}=y_{j_v},\\[6pt]
			\infty, & \text{otherwise}.
		\end{cases}
		\]
	\end{itemize}
\end{remark}

\subsection[A possible link between the sequences in Lemma~\ref{lem:statement}]{A possible link between the sequences in \Cref{lem:statement}}
Although the LDPs in \Cref{teo:low-disorder-1} can in principle be recovered from
\Cref{teo:low-disorder-2} through the contraction principle, the corresponding argument
would be considerably more involved than establishing the LDPs  by a direct and independent approach.

We first assume that $x_1,\ldots,x_m$ are all different from the North Pole $N$.
In this case the sequence
\[
\{\,T_L(\underline{x}) - T_L(N)\,\boldsymbol{1}\,\}_{L\ge 1}
\]
can be viewed as the continuous image of the sequence
\[
\left\{
\left(
T_L(\underline{x}) - \frac{a_{00,L}}{\sqrt{\omega_d}}\,\boldsymbol{1},
\,
T_L(N) - \frac{a_{00,L}}{\sqrt{\omega_d}}
\right)
\right\}_{L\ge 1},
\]
which is a slight modification of the sequence considered in
\Cref{teo:low-disorder-2}.
Define the continuous map
\[
\Psi(z_1,\ldots,z_m,z_{m+1}) := (z_1 - z_{m+1},\,\ldots,\, z_m - z_{m+1}).
\]
Then the contraction principle yields
\begin{multline*}
	I_{1;\mathfrak{g}}(y_1,\ldots,y_m;\underline{x})
	\\
	= \inf
	\left\{
	I_{2;\mathfrak{g}}(z_1,\ldots,z_m,z_{m+1};\,\underline{x},N) :
	\Psi(z_1,\ldots,z_m,z_{m+1}) = (y_1,\ldots,y_m)
	\right\}
	\\[4pt]
	= \inf_{z_{m+1}\in\mathbb{R}}
	\,
	I_{2;\mathfrak{g}}
	(
	y_1 + z_{m+1},\,\ldots,\, y_m + z_{m+1},\, z_{m+1};
	\underline{x},N
	),
\end{multline*}
where $\mathfrak{g}=\mathfrak{L}$ in the low-disorder case and
$\mathfrak{g}=\mathfrak{S}$ in the sparse case.

\medskip

We conclude with a brief discussion of the case in which $x_i = N$ for some
$i\in\{1,\ldots,m\}$.
We assume $m\ge 2$ to avoid trivialities.
In this situation, the argument above can be repeated after removing the $i$-th coordinate,
thus reducing the problem to $m-1$ points.
Indeed, both $I_{1;\mathfrak{L}}(y_1,\ldots,y_m;\underline{x})$ and
$ I_{1;\mathfrak{S}}(y_1,\ldots,y_m;\underline{x})$ are finite only if $y_i = 0$;
if $y_i\neq 0$ then the rate functions are equal to~$\infty$.

\subsection{Some examples}

\begin{example}[Exponential activation function]\label{ex:exp-activation}
	In~\cite{DililloMarinucciSalviVigogna} the authors proved that, for
	\[
	\sigma(x)=e^{-a^2 x^2/2},
	\]
	the corresponding kernel is
	\[
	\kappa(t)
	=\sqrt{\frac{2a^2+1}{(a^2+1)^2-(a^2 t)^2}},
	\qquad t\in[-1,1].
	\]
	Setting $\gamma := \frac{2a^2+1}{a^4}>0$, we may rewrite this as
	\[
	\kappa(t)=\sqrt{\frac{\gamma}{\gamma+1-t^2}},
	\qquad t\in[-1,1].
	\]
	A direct computation shows that $\kappa'(1)=1/\gamma$. Hence we are in the
	low-disorder regime when $\gamma>1$, in the sparse regime when $\gamma=1$, and in
	the high-disorder regime when $\gamma<1$ (see \Cref{ex:high-disorder-example}).
	Notice also that $\mathcal{I}(\kappa)=\{-1,1\}$ in all cases, since $\kappa$ is symmetric
	(see \Cref{lem:set-of-attractive-points-of-iterations}).

	\begin{itemize}
		\item \textbf{$\gamma>1$ (low-disorder).}
		One checks (by induction on $L$) that
		\[
		\kappa_L(t)
		=
		\sqrt{
			\frac{\gamma^{L+1}-\gamma-(\gamma^{L}-\gamma)t^2}{\gamma^{L+1}-1-(\gamma^{L}-1)t^2}}
	\]
	Then, the function defined in \Cref{lem:low-disorder} reduces to,
	\[
	\mathfrak{L}(t)
	:=
	\lim_{L\to\infty} \gamma^L(1-\kappa_L(t))
	=
	\frac{(\gamma-1)(1-t^2)}{2(\gamma-t^2)},
	\qquad t\in[-1,1].
	\]

	\item \textbf{$\gamma=1$ (sparse).}
	Again by induction,
	\[
	\kappa_L(t)
	=
	\sqrt{\frac{L-(L-1)t^2}{L+1-Lt^2}}.
	\]
One can prove that $\kappa_L$ has two derivatives at $t=1$, so
Assumption~\ref{ass:A} holds with $\rho=2$ and $c=2$.
Moreover, \Cref{lem:sparse} yields
\[
\mathfrak{S}(t)
	:=\lim_{L\to\infty} L\,\bigl(1-\kappa_L(t)\bigr)
=
\begin{cases}
	0, & t^2=1\ \text{(i.e., } t\in\mathcal{I}(\kappa)),\\[4pt]
	\frac12, & t^2<1\ \text{(i.e., } t\notin\mathcal{I}(\kappa)).
\end{cases}
\]
\end{itemize}
\end{example}

\begin{example}[ReLU activation function]\label{ex:ReLU}
In~\cite{ChoSaul} the authors proved that for
\[
\sigma(x)=x\vee 0,
\]
the associated kernel is
\[
\kappa(t)
=
\frac{1}{\pi}\Bigl(t(\pi-\arccos t)+\sqrt{1-t^2}\Bigr),
\qquad t\in[-1,1].
\]
In particular, $\kappa'(1)=1$. The kernel is not symmetric and therefore
$\mathcal{I}(\kappa)=\{1\}$ (see \Cref{lem:set-of-attractive-points-of-iterations}).

Unlike the exponential case, it is not easy to obtain a closed expression for $\kappa_L(t)$.
However, one may expand $\kappa$ near $t=1$:
\[
\kappa(t)
=
t
+
\frac{1}{\pi}\!\left(\sqrt{1-t^2}-t\arccos t\right)
=
t+\frac{2\sqrt{2}}{3\pi}(1-t)^{3/2}
+o\!\left((1-t)^{3/2}\right)
\qquad (t\to 1^-).
\]
Indeed,
\[
\lim_{t\to 1^-}
\frac{\sqrt{1-t^2}-t\arccos t}{(1-t)^{3/2}}
=
\frac{2\sqrt{2}}{3},
\]
by l’Hôpital’s rule and the classical limit
\(
\lim_{t\to 1^-}\frac{\arccos t}{\sqrt{1-t}}=\sqrt{2}.
\)
Thus Assumption~\ref{ass:A} holds with $c=\frac{2\sqrt{2}}{3\pi}$ and $\rho=\frac32$.
\end{example}

\section{Functional results}\label{sec:functional}

In this section we investigate to what extent functional LDPs can be established for
one (or both) of the sequences introduced above; see
\Cref{teo:low-disorder-1,teo:low-disorder-2}.
Our aim is to understand whether the finite-dimensional results obtained in the
previous section admit meaningful functional counterparts. We start with the low--disorder regime.

\begin{theorem}\label{prop:low-disorder-functional}
	Assume that $\kappa'(1)<1$.
	Then the sequences $\{(T_L(x)-T_L(N))_{x\in\mathbb{S}^d}\}_{L\ge 1}$ and
	\(
	\left\{
	\left(T_L(x)-\frac{a_{00,L}}{\sqrt{\omega_d}}\right)_{x\in\mathbb{S}^d}
	\right\}_{L\ge 1}
	\)
	satisfy an LDP with the same speed $(\kappa'(1))^{-L}$ and with respective good
	rate functions $I_1$ and $I_2$, given by
    \begin{equation}\label{eq:rf-inf}\small
	I_\ell(f)
	=
	\sup\Bigl\{
	I_{\ell;\mathfrak{L}}\bigl(f(x_1),\ldots,f(x_m);\,x_1,\ldots,x_m\bigr)
	:\ m\ge1,\ x_1,\ldots,x_m\in\mathbb{S}^d,\ x_i\neq x_j\ \text{for } i\neq j
	\Bigr\},
	\end{equation}
	where $\ell=1,2$.
	Here $I_{1;\mathfrak{L}}$ and $I_{2;\mathfrak{L}}$ are the rate functions defined in
	\Cref{teo:low-disorder-1,teo:low-disorder-2}, and $\mathfrak{L}$ is the limiting
	function introduced in \Cref{lem:low-disorder}.
\end{theorem}

It is also well known that $I_1$ and $I_2$ admit a representation in terms of norms
in suitable Reproducing Kernel Hilbert Spaces (see, e.g., \cite{Baldi}).
More precisely, for $\ell\in\{1,2\}$ we have
\[
I_\ell(f)=
\begin{cases}
	\dfrac{\|f\|_{b_\ell}^2}{2}, & \text{if } f\in\mathcal{H}(b_\ell),\\[4pt]
	\infty, & \text{otherwise},
\end{cases}
\]
where $\mathcal{H}(b_\ell)$ is the RKHS generated by the covariance kernel $b_\ell$,
with the associated norm $\|\cdot\|_{b_\ell}$.
The kernels are given by
\[
b_1(z,w)
=
\mathfrak{L}(\langle z,N\rangle)
+\mathfrak{L}(\langle w,N\rangle)
-\mathfrak{L}(\langle z,w\rangle),
\]
and
\[
b_2(z,w)
=
\frac{\displaystyle \int_{-1}^1 \mathfrak{L}(t)(1-t^2)^{\frac{d}{2}-1}\,dt}
{\displaystyle \int_{-1}^1 (1-t^2)^{\frac{d}{2}-1}\,dt}
\;-\;
\mathfrak{L}(\langle z,w\rangle).
\]

We conclude this section with a remark concerning an additional (functional)
weak convergence property of the rescaled fields.
Let $\{(U_L(x))_{x\in\mathbb{S}^d}\}_{L\ge 1}$ denote either of the two sequences
appearing in \Cref{prop:low-disorder-functional}.

\begin{remark}\label{rem:regularity-mathfrak-L}
	The proof of \Cref{prop:low-disorder-functional} (see  Appendix \ref{sect:proofTh4.1}) actually establishes
	the functional convergence in distribution of the sequence of processes
	\[
	\{(\sqrt{v_L} U_L(x))_{x\in\mathbb{S}^d}\}_{L\ge 1}
	\]
	in the space $C(\mathbb{S}^d)$.
	In particular, this sequence converges weakly to a centered continuous Gaussian
	random field.
As a consequence of functional convergence, the limiting covariance function
	must be continuous. Since the finite-dimensional distributions of the limit are
	described by the covariance matrices $B_{1;\mathfrak{L}}(\underline{x})$ and
	$B_{2;\mathfrak{L}}(\underline{x})$ (cf.~\eqref{eq:matrixB1-low-disorder},
	\eqref{eq:matrixB2-low-disorder}), it follows that the associated covariance kernels
	$b_1$ and $b_2$ defined above are continuous functions on $\mathbb{S}^d\times\mathbb{S}^d$.
	In particular, this implies  that the function $\mathfrak{L}$ whose existence is
	established in \Cref{lem:low-disorder} is continuous.
    \end{remark}

Here we present \Cref{prop:sparse-functional}, which shows that no analogue of the
functional results obtained in the low-disorder regime
(cf.~\Cref{prop:low-disorder-functional}) can hold in the sparse case.
In some sense, this is not unexpected.
Indeed, if the family $\{(\sqrt{v_L}U_L(x))_{x\in\mathbb{S}^d}\}_{L\ge 1}$ were to
converge weakly to a centered continuous Gaussian random field with a continuous
covariance function, then the limiting covariance structure appearing in
\Cref{lem:sparse} would necessarily be continuous in the spatial variable.
This would force the function $\mathfrak{S}$ to be continuous.

However, $\mathfrak{S}$ is not continuous; hence neither a functional LDP of the type
established in \Cref{prop:low-disorder-functional} nor the functional weak convergence
result discussed in \Cref{rem:regularity-mathfrak-L} can hold in the sparse regime.

\begin{theorem}\label{prop:sparse-functional}
	Assume that $\kappa'(1)=1$.
	Then  the sequences $\{(T_L(x)-T_L(N))_{x\in\mathbb{S}^d}\}_{L\geq 1}$ and
	\(
	\left\{\left(T_L(x)-\frac{a_{00,L}}{\sqrt{\omega_d}}\right)_{x\in\mathbb{S}^d}\right\}_{L\geq 1}
	\)
do not satisfy an LDP with speed $v_L=L^{1/(\rho-1)}$.
\end{theorem}

\begin{remark}
	\Cref{prop:sparse-functional} shows that the two sequences cannot satisfy an LDP  at the natural speed $v_L = L^{1/(\rho-1)}$.
	One may then ask whether an LDP could hold at a different speed.
At any faster speed, it is easy to see that the two sequences cannot satisfy an
	LDP. At any slower speed, the finite-dimensional distributions do satisfy an LDP,
	but with respect to a degenerate rate function. Consequently, if the two sequences
	were to satisfy an LDP at a slower speed, it would necessarily be governed by a
	degenerate rate function, which is of limited interest in our opinion.
For a detailed discussion of this phenomenon, see Proposition~5.8 in
	\cite{BaldiPacchiarotti}.
\end{remark}

\section{Some remarks on the high-disorder regime \texorpdfstring{$\kappa^\prime(1)>1$}{kappa'(1)>1}}\label{sec:high-disorder}
In this section we take $\kappa^\prime(1)>1$, and we investigate the asymptotic behaviour of the sequence 
$\{\kappa_L(t)\}_{L\geq 1}$ as $L\to\infty$ for every choice of $t\in[-1,1]$. Actually, it is well known that the 
case $t=1$ can be neglected because we have $\kappa_L(1)=1$ for every $L\geq 1$.

It is useful to point out that, when $\kappa^\prime(1)>1$, for $t\in[0,1]$ we have $\kappa(t)\geq 0$, $\kappa$ is 
increasing (and $\kappa(1)=1$); thus there exists $t^*\in[0,1)$ such that $\kappa(t^*)=t^*$. Then we have the 
following situation.

\begin{itemize}
	\item If $t\in[0,1)$, then
	\begin{equation}\label{eq:different-attraction}
		\lim_{L\to\infty}\kappa_L(t)=t^*
	\end{equation}
	because the sequence $\{\kappa_L(t)\}_{L\geq 1}$ is constant for $t=t^*$, decreasing for $t>t^*$ and increasing
	for $t<t^*$ (actually the case $t<t^*$ is empty when $t^*=0$).
	\item If $t\in(-1,0)$, then $\kappa(t)>t$; indeed, by \eqref{eq:kappa-for-sigma}, we have
	$$\kappa(t)-t=\sum_{q\geq 0}\frac{(J_q(\sigma))^2}{q!}t^q
	-\underbrace{\left(\sum_{q\geq 0}\frac{(J_q(\sigma))^2}{q!}\right)}_{=\kappa(1)=1}t
	=\sum_{q\geq 0}\frac{(J_q(\sigma))^2}{q!}(t^q-t)>0,$$
	where the final inequality holds noting that $t^q\geq t$ for every $q\geq 0$, with the strict inequality 
	$t^q>t$ if $q\neq 1$. Therefore the sequence $\{\kappa_L(t)\}_{L\geq 1}$ is increasing and, since we cannot 
	find any fixed point of $\kappa$ in $(-1,0)$, we still get \eqref{eq:different-attraction}.
	\item If $t=-1$, then we can have $\lim_{L\to\infty}\kappa_L(-1)=-1$ if $\kappa(-1)=-1$ or, otherwise,
	$$\lim_{L\to\infty}\kappa_L(-1)=\left\{\begin{array}{ll}
		t^*&\ \mbox{if}\ \kappa(-1)\in(-1,1)\\
		1&\ \mbox{if}\ \kappa(-1)=1.
	\end{array}\right.$$
	\end{itemize}

For completeness we recall that $D_{0,L}\to\widehat{d}_0<1$ (as $L\to\infty$); see
Theorem 3.3 in \cite{DililloMarinucciSalviVigogna}. Thus the asymptotic analysis of $\{\kappa_L(t)\}_{L\geq 1}$ 
presented above allows to say that the two limit at zero in Lemma \ref{lem:statement} fail. Actually this was
already stated in Theorem 3.3 in \cite{DililloMarinucciSalviVigogna} for $\{T_L(x)-T_L(N)\}_{L\geq 1}$.

One could try to investigate the functional convergence of the sequences in Lemma \ref{lem:statement}.
In particular one could start with a large class of finite dimensional distributions, for instance for 
$x_1,\ldots,x_m\in\mathbb{S}^d$ such that $|\langle x_i,x_j\rangle|<1$ for $i\neq j$.
\begin{enumerate}
	\item For simplicity here we also assume that $|\langle x_i,N\rangle|<1$ for every  $i\in\{1,\ldots,m\}$.
	Then we have
	$$\lim_{L\to\infty}\mathrm{\mathrm{Cov}}(T_L(x_i)-T_L(N),T_L(x_j)-T_L(N))=\left\{\begin{array}{ll}
		2(1-t^*)&\ \mbox{if}\ i=j\\
		1-t^*&\ \mbox{if}\ i\neq j.
	\end{array}\right.$$
	\item We have
	$$\lim_{L\to\infty}\mathrm{\mathrm{Cov}}\left(T_L(x_i)-\frac{a_{00,L}}{\sqrt{\omega_d}},T_L(x_j)-\frac{a_{00,L}}{\sqrt{\omega_d}}\right)=\left\{\begin{array}{ll}
		1-\widehat{d}_0&\ \mbox{if}\ i=j\\
		t^*-\widehat{d}_0&\ \mbox{if}\ i\neq j.
	\end{array}\right.$$
\end{enumerate}
Then, in both cases, we cannot have the functional weak convergence.  Indeed,
if one of the two sequences does converge weakly to a centered continuous Gaussian random field, with a 
suitable (continuous) covariance function, this would be in contradiction with the computations above for the covariance 
functions for the finite dimensional distributions.

We conclude with an example.

\begin{example}[Example \ref{ex:exp-activation} with $\gamma<1$]\label{ex:high-disorder-example}
	One can check (by induction on $L$) that
	$$\kappa_L(t)=\sqrt{\frac{\gamma^{L+1}-\gamma-(\gamma^L-\gamma)t^2}{\gamma^{L+1}-1-(\gamma^L-1)t^2}}.$$
	Then it is easy to check that $t^*$ above
	coincides with $\sqrt{\gamma}$, and therefore
	$$\lim_{L\to\infty}\kappa_L(t)=\left\{\begin{array}{ll}
		\sqrt{\gamma}&\ \mbox{if}\ t\in(-1,1)\\
		1&\ \mbox{if}\ |t|=1.
	\end{array}\right.$$
\end{example}

\appendix
\section{Proofs}\label{app}
In this section we provide the proofs of all the statements given above.

\subsection{Proof of Lemma \ref{lem:set-of-attractive-points-of-iterations}}
Using the triangle inequality in~\eqref{eq:kappa-for-sigma}, we obtain for every $t\in (-1,1)$
\begin{equation}\label{eq:open-interval-property}
	\Big|\sum_{q\geq 0}  \frac{ (J_q(\sigma))^2}{q!} t^q\Big|
	\;\leq\; \sum_{q\geq 0} \frac{ (J_q(\sigma))^2}{q!} |t|^q
	\;<\; \sum_{q\geq 0}  \frac{ (J_q(\sigma))^2}{q!} = \kappa(1)=1,
\end{equation}
where the strict inequality holds because $|t|<1$ and there exists at least one index $q\ge 1$ such that $J_q(\sigma)\neq 0$.

From~\eqref{eq:open-interval-property}, we have  $\mathcal I(\kappa)\subseteq \{ -1,1 \}$, and obviously  $1\in\mathcal{I}(\kappa)$.  We note that $\kappa(-1)  = 1$ if and only if $\kappa$ is even function. Indeed, from
\begin{align*}
&1=\kappa(1)=\sum_{q\geq 0} \frac{ (J_{2q}(\sigma))^2}{(2q)!} +\sum_{q\geq 0} \frac{ (J_{2q+1}(\sigma))^2}{(2q+1)!}\\
&1=\kappa(-1)=\sum_{q\geq 0} \frac{ (J_{2q}(\sigma))^2}{(2q)!}-\sum_{q\geq 0} \frac{ (J_{2q+1}(\sigma))^2}{(2q+1)!}
\end{align*}
we have  $J_{2q+1}(\sigma)=0$ for every $q\geq 0$.

\subsection{Proof of Lemma \ref{lem:low-disorder}}
Before proving the lemma, we first show that if $\kappa'(1)<1$, then
\begin{equation}\label{eq:mon-derivata}
|\kappa'(t)| < \kappa'(1) \qquad \text{for every } t\in (-1,1).
\end{equation}
From \eqref{eq:kappa-for-sigma}, the derivative $\kappa'(t)$ for $t\in(-1,1)$ can be obtained
term-by-term from the power series.
Moreover, by the monotone convergence theorem we have
\begin{equation}\label{eq:monotone-convergence-consequence}
	\kappa'(1) = \lim_{t\to 1^-} \kappa'(t).
\end{equation}
The bound on $|\kappa'(t)|$ follows in exactly the same way as in
\eqref{eq:open-interval-property}.

Let us prove the first part of the lemma, namely that the sequence
$$ \beta_L(t) := \frac{1-\kappa_L(t)}{(\kappa'(1)^L)}
$$
has limit in $[0,\infty)$.
Since $\kappa_L(1)=1$ for every $L\geq 1$, we trivially obtain $\mathfrak{L}(1)=0$.
Let $t\neq 1$.
Because $\kappa\in C^1([-1,1])$, there exists $\xi_{L,t}\in (t,1)\subset (-1,1)$ such that
\[
0\leq 1-\kappa_L(t)
=1-\kappa(\kappa_{L-1}(t))
=\kappa'(\xi_{L,t})\bigl(1-\kappa_{L-1}(t)\bigr)
\le \kappa'(1)\bigl(1-\kappa_{L-1}(t)\bigr),
\]
where the inequality follows from~\eqref{eq:mon-derivata}.
Hence $0\leq \beta_L(t)\leq \beta_{L-1}(t)$.

We therefore obtain a monotone and bounded sequence, which guarantees the existence
(and finiteness) of the limit $\mathfrak{L}(t)$.

The proof of the second part of the lemma, i.e. $\mathfrak L(t) =0$ if and only if $t\in \mathcal I(\kappa)$, is inspired by Theorem 1.3.1 in
\cite{KuczmaChoczewskiGer}. Trivially for $t\in \mathcal{I}(\kappa)$ we have $\mathfrak L(t)=0$. In what follows, we assume that $t\not\in \mathcal I(\kappa)$.
We note that for $L\geq 1$
\begin{equation}\label{eq:numero}\beta_L(t)=(1-t)\prod_{\ell=0}^{L-1}\left(1+\frac{p_\ell(t)}{\kappa^\prime(1)}
\right) = (1-t)\exp\left(\sum_{\ell=0}^{L-1}\log\left(1+\frac{p_\ell(t)}{\kappa^\prime(1)}\right)
\right)
\end{equation}
where
$$p_\ell(t):=\frac{1-\kappa_{\ell+1}(t)}{1-\kappa_\ell(t)}-\kappa^\prime(1)\quad\mbox{for}\ \ell\geq 0,$$
with $\kappa_0(t)=t$ and we note that the second equality in~\eqref{eq:numero} is well defined since $p_\ell(t)\geq-\kappa^\prime(1)$.
We conclude the proof if we show that the series
$\sum_{\ell=0}^\infty\log\left(1+\frac{p_\ell(t)}{\kappa^\prime(1)}\right)$ is convergent.

By Lagrange's Theorem, there exists $\xi_{\ell,t} \in (\kappa_\ell(t),1)$ such that
\[
p_\ell(t)
= \frac{1 - \kappa_{\ell+1}(t)}{1 - \kappa_\ell(t)} - \kappa'(1)
= \frac{1 - \kappa(\kappa_\ell(t))}{1 - \kappa_\ell(t)} - \kappa'(1) =  \kappa'(\xi_{\ell,t}) - \kappa'(1).
\]

Since, as $\ell\to + \infty$,  $\kappa_\ell(t) \to 1$, we have $\xi_{\ell,t} \to 1$, hence $\kappa'(\xi_{\ell,t}) \to \kappa'(1)$.
Thus, for every $\eta>0$, there exists $\bar{\ell}$ such that, for every $\ell>\bar{\ell}$,
$$\left|\log\left(1+\frac{p_\ell(t)}{\kappa^\prime(1)}\right)\right|\leq(1+\eta)\frac{|p_\ell(t)|}{\kappa^\prime(1)},$$
and it suffices to show that $\sum_{L=1}^\infty|p_L(t)|<\infty$. Using Assumption~\ref{ass:A}, for $\ell$ large enough, we have
\[
p_\ell(t)
=\frac{1-\kappa(\kappa_\ell(t))}{1-\kappa_\ell(t)}-\kappa^\prime(1)\\
=-c(1-\kappa_\ell(t))^{\rho-1}+o((1-\kappa_\ell(t))^{\rho-1});
\]
then
$$|p_\ell(t)|\leq 2|c|(1-\kappa_\ell(t))^{\rho-1}.$$
Finally,  we have
$$(1-\kappa_\ell(t))^{\rho-1}\leq\{\kappa^\prime(1)(1-\kappa_{\ell-1}(t))\}^{\rho-1}
\leq\cdots\leq(\kappa^\prime(1)^{\rho-1})^{\ell-\bar{\ell}}(1-\kappa_{\bar{\ell}}(t))^{\rho-1};$$
in conclusion, at least for $\ell$ large enough, $|p_\ell(t)|$ is majorized by the term of a convergent geometric series
(multiplied by a constant), and this completes the proof.

\subsection{Proof of Lemma \ref{lem:sparse}}
Firstly, note that, by $\kappa^\prime(1)=1$, \eqref{eq:sort-of-Taylor-formula} yields
\begin{equation}\label{eq:sort-of-Taylor-formula-consequence}
	\kappa(t)=t+c(1-t)^\rho+o((1-t)^\rho)\quad\mbox{as $t\to 1^-$}
\end{equation}
for $c\neq 0$ and $\rho>1$. Let us prove $c>0$. To do this, we prove that  $\kappa(t)>t$ for $|t|<1$.
We consider the function
$$g(t):=\kappa(t)-t;$$
then we have $g^\prime(t)=\kappa^\prime(t)-1$. We note that  \eqref{eq:monotone-convergence-consequence} holds also if $\kappa'(1)=1$ and hence
$g^\prime(t)\leq\kappa^\prime(1)-1\leq 0$, and the inequality is strict for $|t|<1$; so $g$ is a strictly
decreasing function. Thus, for $|t|<1$, we have $g(t)>g(1)$, i.e., $\kappa(t)-t>\kappa(1)-1=0$, or equivalently
$\kappa(t)>t$, 
then $c>0$.

Now we discuss the limit in the statement of the lemma, which trivially holds if
$t\in\mathcal{I}(\kappa)$; so, from now on, in this proof, we take $t\notin\mathcal{I}(\kappa)$.
By \eqref{eq:sort-of-Taylor-formula-consequence} (and since $\kappa_L(t)\to 1$ as $L\to\infty$) we have
\begin{multline*}
	1-\kappa_L(t)=1-\kappa(\kappa_{L-1}(t))=1-\{\kappa_{L-1}(t)+c(1-\kappa_{L-1}(t))^\rho+o((1-\kappa_{L-1}(t))^\rho)\}\\
	=1-\kappa_{L-1}(t)-c(1-\kappa_{L-1}(t))^\rho+o((1-\kappa_{L-1}(t))^\rho)\\
	=(1-\kappa_{L-1}(t))\{1-c(1-\kappa_{L-1}(t))^{\rho-1}+o((1-\kappa_{L-1}(t))^{\rho-1})\}
\end{multline*}
and therefore (we recall that $(1+z)^{-\rho+1}\sim 1+(-\rho+1)z$ as $z\to 0$)
$$(1-\kappa_L(t))^{-\rho+1}=(1-\kappa_{L-1}(t))^{-\rho+1}+c(\rho-1)+o(1).$$
Then
$$\lim_{L\to\infty}(1-\kappa_L(t))^{-\rho+1}-(1-\kappa_{L-1}(t))^{-\rho+1}=c(\rho-1)$$
and, by Stolz-Cesaro Theorem,
$$\lim_{L\to\infty}\frac{(1-\kappa_L(t))^{-\rho+1}}{L}=c(\rho-1),$$
or equivalently
$$\lim_{L\to\infty}\frac{1-\kappa_L(t)}{L^{-1/(\rho-1)}}=(c(\rho-1))^{-1/(\rho-1)}.$$
This completes the proof.

\subsection[Proof of Theorem~\ref{teo:low-disorder-1}]{Proof of \Cref{teo:low-disorder-1}}

In this proof we set $$U_L(\underline{x})=T_L(\underline x )-T_L(N)\boldsymbol{1}.$$ Hence, since $U_L(\underline x)$ is a zero mean Gaussian random with covariance
$$  \mathrm {Cov}(U_L(x_i) , U_L(x_j)) = \kappa_L(\langle x_i, x_j\rangle ) - \kappa_L(\langle x_i, N\rangle ) - \kappa_L(\langle x_j, N\rangle ) +1, $$
for every $\underline{\theta}\in\mathbb{R}^m$ we have
$$\frac{1}{v_L}\log\mathbb{E}\left[e^{v_L\langle\underline{\theta}, U_L(\underline{x}) \rangle}\right]\\
=\frac{v_L}{2}\sum_{i,j=1}^m\{\kappa_L(\langle x_i,x_j\rangle)-\kappa_L(\langle x_i,N\rangle)
-\kappa_L(\langle x_j,N\rangle)+1\} \theta_{i} \theta_j.$$
Now, using Lemma \ref{lem:low-disorder} for the low-disorder cases and Lemma \ref{lem:sparse} for the sparse case, we have
$$ \lim_{L\to + \infty } v_L \big( \kappa_L(\langle x_i,x_j\rangle)-\kappa_L(\langle x_i,N\rangle)
-\kappa_L(\langle x_j,N\rangle)+1\big)  = \mathfrak g(\langle x_i, N\rangle)  + \mathfrak{g}(\langle x_j, N\rangle) - \mathfrak{g}(\langle x_i, x_j\rangle)$$
where $v_L$ and $\mathfrak g$ are as in the statement of the theorem. Thus,
$$\lim_{L\to\infty}\frac{1}{v_L}\log\mathbb{E}\left[e^{v_L\langle\underline{\theta},U_L(\underline{x})\rangle}\right]
=\frac{1}{2}\langle\underline{\theta},B_{1;\mathfrak g}(\underline{x})\underline{\theta}\rangle,$$
where $B_{1;\mathfrak g}(\underline{x})$ is the square matrix in \eqref{eq:matrixB1-low-disorder}.
Then the desired LDP holds by applying the G\"artner Ellis Theorem.

\subsection[Proof of Theorem~\ref{teo:low-disorder-2}]{Proof of \Cref{teo:low-disorder-2}}
In the proof of the theorem we shall need to interchange a limit and an integral.
This step is justified by the following lemma, whose proof is postponed to the end of this section.

\begin{lemma}\label{lem:aux}Let us assume $\kappa'(1)\leq 1$. Then
	$$\lim_{L\to + \infty }\frac{\int_{-1}^1 v_L(1-\kappa_L(t))(1-t^2)^{\frac{d}{2}-1}dt}{\int_{-1}^1(1-t^2)^{\frac{d}{2}-1}dt}=\frac{\int_{-1}^1\mathfrak{g}(t)(1-t^2)^{\frac{d}{2}-1}dt}{\int_{-1}^1(1-t^2)^{\frac{d}{2}-1}dt}$$
	where $v_L$ and $\mathfrak{g}$ are as in~\Cref{teo:low-disorder-1}.
\end{lemma}

In this proof we set $$U_L(\underline{x})=T_L(\underline x )-\frac{a_{00,L}}{\sqrt{\omega_d}}\boldsymbol 1.$$ Since the sequence $\{a_{\ell,m}\}_{(\ell,m)\in \mathcal A}$ are independent centered Gaussian random variable with variance $\mathbb E[a_{\ell m,L}^2] = C_{\ell,L}$, we have  
$$ \mathrm{Cov}\left(T_L(x_i), \frac{a_{00,L}}{\sqrt{\omega_d} }\right) = D_{0,L} $$
and hence
$$  \mathrm {Cov}(U_L(x_i) , U_L(x_j)) = \kappa_L(\langle x_i, x_j\rangle ) - D_{0,L}.$$
Thus, for every $\underline{\theta}=(\theta_1,\ldots,\theta_m)\in\mathbb{R}^m$ we have
\begin{align*}\frac{1}{v_L}\log\mathbb{E}\left[e^{v_L\langle\underline{\theta},U_L(\underline{x})\rangle}\right]
&=\frac{v_L}{2}\sum_{i,j=1}^m\{\kappa_L(\langle x_i,x_j\rangle)-D_{0,L}\}\theta_i \theta_j \,.
\end{align*}
Now, 
$$ v_L \big( \kappa_L(\langle x_i, x_j\rangle) -D_{0,L}\big)  =  v_L (1 - D_{0,L}) - v_L ( 1-\kappa_L(\langle x_i, x_j \rangle) )\,.$$ 
Using the orthogonality of the Gegenbauer polynomials, we have
$$ D_{0,L} = \frac{\int_{-1}^1 \kappa_L(t) (1-t^2)^{d/2-1} \mathrm dt }{\int_{-1}^1(1-t^2)^{d/2-1} \mathrm d t } $$
and hence
$$ v_L( 1- D_{0,L}) =\frac{\int_{-1}^1 (1-\kappa_L(t)) (1-t^2)^{d/2-1} \mathrm dt }{\int_{-1}^1(1-t^2)^{d/2-1} \mathrm d t }  \,. $$

Thus, combining \Cref{lem:aux} with \Cref{lem:low-disorder}, for the low-disorder case, and \Cref{lem:sparse}, for the sparse, we obtain
$$ \lim_{L\to + \infty}  v_L \big( \kappa_L(\langle x_i, x_j\rangle -D_{0,L}\big)  =  	\frac{\int_{-1}^1\mathfrak{g}(t)(1-t^2)^{\frac{d}{2}-1}dt}
	{\int_{-1}^1(1-t^2)^{\frac{d}{2}-1}dt}-\mathfrak{g}(\langle x_i,x_j\rangle). $$
Then the desired LDP holds by applying the
G\"artner Ellis Theorem.

\medskip

We split the proof of \Cref{lem:aux} into two cases.
\begin{proof}[Proof of~\Cref{lem:aux} for the low-disorder case] Using~\Cref{lem:low-disorder} we have
$$ \beta_L(t) = \frac{ 1- \kappa_L(t)}{(\kappa'(1))^L}  \to \mathfrak{L}(t) \ . $$
Since $|\kappa'(t)|< 1$ for every $t\in (-1,1)$, using $L$ times the mean value theorem we have
\begin{equation}\label{bound}0\leq \beta_L(t)\leq 1-t  \end{equation}
thus $\beta_L\leq 2 $ and the claim follows from the dominated convergence.
\end{proof}

\begin{proof}[Proof of~\Cref{lem:aux} for the sparse case] As in the previous case, using~\Cref{lem:sparse}
	$$ \eta_L(t): = L^{1/(\rho-1)} ( 1- \kappa_L(t))
\to \mathfrak S(t). $$
In order to apply the dominated convergence, we prove that there exists a constant $M>0$ such that $ 0 \leq \eta_L(t)\leq M $. For every $L\geq 1$ let
$$R_L:=\{t\in[-1,1]:\kappa_L(t)>0\}.$$
By construction, $\{ R_L\}_{L\geq 1}$ is an open cover of $[-1,1]$, indeed  $\kappa_L(t) \to 1$ for every $t\in [-1,1]$.
Then, since $[-1,1]$ is a compact set,
there exists a finite subcover of $[-1,1]$; thus, for some $\bar{L}$, we have
$$[-1,1]=\bigcup_{i=1}^{\bar L} R_i $$ and since
for $t\in(0,1]$, we have $\kappa(t)>0$ we have $R_L\subset R_{L+1}$. Thus $[-1,1] = R_{\bar L}$.
Now we take $L>\bar{L}$ and, for every $t\in[-1,1]$, we get
$$\kappa_L(t)=\kappa_{L-\bar{L}}(\kappa_{\bar{L}}(t))\geq\kappa_{L-\bar{L}}(0),$$
whence we obtain
$$0\leq L^{1/(\rho-1)}(1-\kappa_L(t))\leq\frac{L^{1/(\rho-1)}}{(L-\bar{L})^{1/(\rho-1)}}
(L-\bar{L})^{1/(\rho-1)}(1-\kappa_{L-\bar{L}}(0)),$$
and the upper bound converges to $\mathfrak{S}(0)$ (as $L\to\infty$) by Lemma \ref{lem:sparse}. In conclusion the
proof is complete because a convergent sequence is bounded.
 	\end{proof}

\subsection{Proof of \Cref{prop:low-disorder-functional}} \label{sect:proofTh4.1}
Throughout this proof we denote by $\{(U_L(x))_{x\in S^d}\}_{L\ge1}$ one of the two sequences in the statement of the theorem.
 In both cases, the finite-dimensional distributions satisfy an LDP 
(see \Cref{teo:low-disorder-1,teo:low-disorder-2}). 
By the Dawson G\"artner Theorem, 
any LDP for the sequence $\{(U_L(x))_{x\in S^d}\}_{L\ge1}$ in $C(\mathbb S^d)$, if it exists, 
must have rate function given by \eqref{eq:rf-inf}. 
By classical results on LDP, to prove the theorem it suffices to show that 
the sequence $\{(U_L(x))_{x\in\mathbb S^d}\}_{L\ge 1}$ is exponentially tight in $C(\mathbb S^d)$. 
That is, for every $M>0$ there exists a compact set
$K_M \subset C(\mathbb S^d)$
such that
$$
\limsup_{L\to\infty}\frac{1}{v_L}\log \mathbb P((U_L(x))_{x\in\mathbb S^d}\notin K_M)\le -M.
$$

Using Theorem 2.4 in \cite{Baldi}, to prove that the sequence $\{(U_L(x))_{x\in\mathbb{S}^d}\}_{L\ge 1}$ is exponentially tight at the speed $v_L$, 
it is sufficient to show that $\{(\sqrt{v_L} U_L(x))_{x\in\mathbb{S}^d}\}_{L\ge 1}$ is tight in $C(\mathbb{S}^d)$; 
the interested reader can see \cite{Baldi} for more details on the relationship between tightness and exponential tightness. 

From the weak convergence of the finite-dimensional distributions of the sequence
$$
\{(\sqrt{v_L}\,U_L(x))_{x\in\mathbb{S}^d}\}_{L\ge 1},
$$
if one can show that, for some constants $a,b>0$ and for every $x,y\in \mathbb S^d$ and every $L\ge 1$,
$$
\mathbb{E}[|\sqrt{v_L} U_L(x)-\sqrt{v_L} U_L(y)|^a] \le c \, \|x-y\|^{b+d},
$$
then, by Theorem 23.7 in \cite{Kallenberg}, the sequence 
$
\{(\sqrt{v_L}\,U_L(x))_{x\in\mathbb{S}^d}\}_{L\ge 1}
$
is tight in $C(\mathbb S^d)$ and its limit has continuous sample paths.

Let us prove that the previous inequality holds with $a = 2r$ and $b = 2r-d$ for $r>0$ large enough. 
We note that, in both cases,
$$
U_L(x) - U_L(y) = T_L(x) - T_L(y)  \stackrel{d}{=} 
 \sqrt{2(1-\kappa_L(\langle x, y\rangle ))} \, Z,
$$
where $\stackrel{d}{=} $ denotes equality in law and $Z$ is a standard Normal random variable. Thus,
$$
\mathbb{E}[|\sqrt{v_L}U_L(x)-\sqrt{v_L}U_L(y)|^{2r}] = (2v_L(1-\kappa_L(\langle x,y\rangle)))^r \, \mathbb{E}[|Z|^{2r}].
$$

Now, since $v_L=(\kappa^\prime(1))^{-L}$, recalling that 
$v_L (1- \kappa_L(\langle x,y\rangle)) = \beta_L(\langle x,y\rangle)$ 
and using \eqref{bound}, we obtain
$$
2v_L(1-\kappa_L(\langle x,y\rangle)) \le 2 (1-\langle x,y\rangle) = \|x-y\|^2,
$$
where we have used that $\|x\| = \|y\| = 1$. In conclusion, we have
$$
\mathbb{E}[|\sqrt{v_L}U_L(x)-\sqrt{v_L}U_L(y)|^{2r}] \le \mathbb{E}[|Z|^{2r}] \, \|x-y\|^{2r},
$$
which yields the desired estimate.

\subsection{Proof of \texorpdfstring{\Cref{prop:sparse-functional}}{prop.~\ref{prop:sparse-functional}}}
We reason by contradiction for both sequences. If functional LDPs hold with speed $v_L=L^{1/(\rho-1)}$,
 in analogy with 
the LDPs stated in  \Cref{prop:low-disorder-functional} (obviously here $\mathfrak{S}$ is the function in Lemma
\ref{lem:sparse}), one should have the following statements.
	\begin{enumerate}
		\item $\{(T_L(x)-T_L(N))_{x\in\mathbb{S}^d}\}_{L\geq 1}$ satisfies the LDP with speed 
		$v_L=L^{1/(\rho-1)}$ and  rate function $I_1$ defined by
		$$I_1(f)=\sup\{I_{1;\mathfrak{S}}(f(x_1),\ldots,f(x_m);x_1,\ldots,x_m):m\geq 1,\
		x_1,\ldots,x_m\in\mathbb{S}^d,\ x_i\neq x_j\ \mbox{for}\ i\neq j\},$$
		\item $\left\{\left(T_L(x)-\frac{a_{00,L}}{\sqrt{\omega_d}}\right)_{x\in\mathbb{S}^d}
		\right\}_{L\geq 1}$ satisfies the LDP with speed $v_L=L^{1/(\rho-1)}$ and  rate 
		function $I_2$ defined by
		$$I_2(f)=\sup\{I_{2;\mathfrak{S}}(f(x_1),\ldots,f(x_m);x_1,\ldots,x_m):m\geq 1,\
		x_1,\ldots,x_m\in\mathbb{S}^d,\ x_i\neq x_j\ \mbox{for}\ i\neq j\}.$$
	\end{enumerate}
    We shall see that $I_i$ are  good rate function and $I_i(f) = +\infty$ for a class of functions that is too large (with a distinction between the cases $i = 1$ and $i = 2$; see below). Therefore, if we try to obtain the finite-dimensional rate function by applying the contraction principle, we would expect
\begin{equation}\label{eq:marginalization-on-fd}
I_{i;\mathfrak{S}}(y_1,\ldots,y_m;x_1,\ldots,x_m) = \inf\{ I_i(f) : f(x_1)=y_1,\ldots,f(x_m)=y_m \}
\end{equation}
for every choice of $x_1,\ldots,x_m \in \mathbb{S}^d$ and $y_1,\ldots,y_m \in \mathbb{R}$. However, this leads to a contradiction, since the infimum above is often $+\infty$, while the actual finite-dimensional rate functions are finite.

\begin{enumerate}
		\item Let $f \in C(\mathbb{S}^d)$ be a non-constant function. Then there exists a closed set $U \subseteq \mathbb{S}^d$ such that $f$ is non-constant on $U$ and $U$ does not contain antipodal points. If we set 
$m_f := \min\{ f(x) : x \in U \}$ and $M_f := \max\{ f(x) : x \in U \}$, we have $m_f < M_f$.
Moreover, it is easy to show that there exists $\varepsilon > 0$ such that, for every $m$, we can find distinct points $x_1, \ldots, x_m \in \mathbb{S}^d$ that
\[
\frac{1}{m} \sum_{i=1}^m f(x_i) = \frac{m_f + M_f}{2}
\]
and
\[
\left| f(x_1) - \frac{m_f + M_f}{2} \right|, \ldots, \left| f(x_m) - \frac{m_f + M_f}{2} \right| > \varepsilon.
\]
(This is easy to do by taking $m$ even, but this is not necessary; moreover, this assumption is not restrictive for the proof.)
 Then, taking into account the discussion in Remark \ref{rem:rf-sparse-1} for the rate function in 
		\Cref{teo:low-disorder-2}, we have
		\begin{eqnarray*}
			I_1(f)&\geq&\frac{1}{2(m+1)h_{c,\rho}}\left\{m\sum_{i=1}^m(f(x_i))^2-2\sum_{1\leq i<j\leq m}f(x_i)f(x_j)\right\}\\
			&=&\frac{1}{2(m+1)h_{c,\rho}}\left\{(m+1)\sum_{i=1}^m(f(x_i))^2-\left(\sum_{i=1}^mf(x_i)\right)^2\right\}\\
			&=&\frac{1}{2h_{c,\rho}}\left\{\sum_{i=1}^m(f(x_i))^2-\frac{1}{m+1}\left(\sum_{i=1}^mf(x_i)\right)^2\right\}\\
			&\geq&\frac{1}{2h_{c,\rho}}\left\{\sum_{i=1}^m(f(x_i))^2-\frac{1}{m}\left(\sum_{i=1}^mf(x_i)\right)^2\right\},
		\end{eqnarray*}
		and therefore
		$$I_1(f)\geq\frac{1}{2h_{c,\rho}}\sum_{i=1}^m\left(f(x_i)-\frac{1}{m}\sum_{j=1}^mf(x_j)\right)^2
		\geq\frac{m\varepsilon^2}{2h_{c,\rho}}.$$
		Given that the above construction is possible for every $m \in \mathbb{N}$, it follows that $I_1(f) = +\infty$.
Moreover, to apply the contraction principle, we need to verify that $I_1$ is a good rate function. In particular, we must show that every level set
\[
\{ f : I_1(f) \leq \eta \}, \quad \eta \geq 0,
\]
which is reduced to constant functions, is a compact set. We assume that there exists $c\in\mathbb{R}$ 
		such that $f(x)=c$ for all $x\in\mathbb{S}^d$. Then, if we restrict the attention to distinct points 
		$x_1,\ldots,x_m\in\mathbb{S}^d$ such that $\langle x_i,x_j\rangle\notin\mathcal{I}(\kappa)$ for every 
		$i,j\in\{1,\ldots,m\}$, we have
		$$I_1(f)\geq\sup_{m\geq 1}\frac{1}{2(m+1)h_{c,\rho}}\left(m^2c^2-2\frac{m(m-1)}{2}c^2\right)
		=\sup_{m\geq 1}\frac{mc^2}{2(m+1)h_{c,\rho}}=\frac{c^2}{2h_{c,\rho}}.$$
		Thus the level set, which we know is closed, is also bounded because it is a subset of
		$\left\{c\in\mathbb{R}:\frac{c^2}{2h_{c,\rho}}\leq\eta\right\}$ (which is bounded); then we can say that
		$I_1$ is a good rate function. In conclusion \eqref{eq:marginalization-on-fd} fails by taking 
		$y_1,\ldots,y_m$ not all coincident.
		\item Take $f \in C(\mathbb{S}^d)$ such that $f(x) \neq 0$ for some $x \in \mathbb{S}^d$. Then there exists $\varepsilon > 0$ such that, for every $m$, we can find $x_1,\ldots,x_m \in \mathbb{S}^d$ with 
\[
|f(x_1)|,\ldots,|f(x_m)| > \varepsilon.
\]
By taking into account the discussion in Remark~\ref{rem:rf-sparse-2} for the rate function in Theorem~\ref{teo:low-disorder-2} (restricting to distinct points $x_1,\ldots,x_m$ such that 
$\langle x_i, x_j \rangle \notin \mathcal{I}(\kappa)$ for every $i,j$; this causes no issues in view of what we will prove), we have
\[
I_2(f) \geq \frac{1}{2h_{c,\rho}} \sum_{i=1}^m (f(x_i))^2 \geq \frac{m \varepsilon^2}{2h_{c,\rho}}.
\]
Since this construction is possible for every $m \in \mathbb{N}$, it follows that $I_2(f) = +\infty$.  
Moreover, if $f(x) = 0$ for all $x \in \mathbb{S}^d$, then $I_2(f) = 0$ (this can be easily checked), and therefore
\[
I_2(f) =
\begin{cases}
0 & \text{if $f(x) = 0$ for all $x \in \mathbb{S}^d$,}\\
+\infty & \text{otherwise.}
\end{cases}
\]
Thus, for every $\eta \geq 0$, the level set $\{ f : I_2(f) \leq \eta \}$ is compact because it reduces to a single point (the null function). Therefore, $I_2$ is a good rate function. In conclusion, \eqref{eq:marginalization-on-fd} fails when taking $y_1,\ldots,y_m$ not all equal to zero.
        
	\end{enumerate}

\bibliographystyle{amsplain} 
\bibliography{references}
\end{document}